\documentclass[hyperref]{ctexart}
\usepackage{geometry} 
\usepackage{helvet}
\usepackage{amsmath, amsfonts, amssymb}
\usepackage[english]{babel}

\usepackage{url}      
\usepackage{bm}      
\usepackage{multirow}
\usepackage{booktabs}
\usepackage{cite}
\usepackage{lipsum}
\usepackage{graphicx}
\usepackage{graphics}
\usepackage{subfigure}
\usepackage{epstopdf}
\usepackage{float}
\usepackage{indentfirst}
\usepackage{algorithm}
\usepackage{algorithmicx}  
\usepackage{algpseudocode} 

\usepackage{fancyhdr} 
\usepackage{hyperref} 
\hypersetup{colorlinks, bookmarks, unicode} 

\numberwithin{equation}{section}

\newtheorem{Theorem}{Theorem}[section]

\newtheorem{Property}{Property}[section]
\newtheorem{Proof}{Proof}[section]
\newtheorem{Remark}{Remark}[section]
\newtheorem{Monotonicity Criterion}{Monotonicity Criterion}
\newtheorem{Corollary}{Corollary}[section]
\newtheorem{Definition}{Definition}[section]
\newtheorem{Proposition}{Proposition}[section]
\newtheorem{Assumption}{Assumption}[section]
\newtheorem{Note}{Note}[section]

\title{\textbf{Riemann problem of Euler equations with singular sources}}
\author{\sffamily Changsheng Yu, \sffamily Tiegang Liu, \sffamily Chengliang Feng}

\begin{document}
	\maketitle
	
\begin{abstract}
	\begin{abstract}
		This paper is concerned with the Riemann problem of one-dimensional Euler equations with a singular source.
		The exact solution of this Riemann problem contains a stationary discontinuity induced by the singular source, which is different from all the simple waves in the Riemann solution of classical Euler equations.
		We propose an eigenvalue-based monotonicity criterion to select the physical curve of this stationary discontinuity.
		By including this stationary discontinuity as an elementary wave, the structure of Riemann solution becomes diverse, e.g. the number of waves is not fixed and interactions between two waves become possible.
		Under the double CRP framework, we prove all possible structures of the Riemann solution.
	\end{abstract}
\end{abstract}
	
\noindent{\bf Keywords: } hyperbolic conservation law, Euler equation, Riemann problem, singular source, entropy condition

\noindent{\bf AMS subject classifications: }35L67,35Q31,76N30

\section{Introduction}\label{introduction}
In governing equations for fluids, the effects other than convection and diffusion are usually described by source terms. 
Sources can exhibit singularities in some complex problems, such as surface tension \cite{fechter2018approximate,houim2013ghost}, evaporation\cite{das2020sharp,das2020simulation,das2021sharp}, phase transition\cite{hitz2020comparison,lee2017sharp,long2021fully,paula2019analysis}, condensation\cite{cheng2010condensation}, geometric constraints\cite{coquel2014robust,lefloch2011godunov,li2021well}.
The singular source induces jumps in the state of flow and leads to new difficulties in the analysis of fluid motion.

This paper focuses on investigating the effect of singular sources on flow field for inviscid compressible fluids. The governing equations are
\begin{equation}\label{equ: governing equation}
	\partial _tU+\partial _xF(U)=\delta(x)S,
\end{equation}
where
\[
U=
\begin{pmatrix}
	\rho\\\rho u\\E
\end{pmatrix},
F=
\begin{pmatrix}
	\rho u\\\rho u^2+p\\(E+p)u
\end{pmatrix},
\delta(x)=\begin{cases}
	0,\ x\ne 0\\
	\infty,\ x=0
\end{cases}
\]
Here, $\rho$, $p$ and $E$ denote the density, pressure and total energy, respectively. $u$ is the velocity. $S$ is the vector of souce term. $\delta(x)$ is the Dirac delta function, which means that the source term is only distributed at the origin. 

Some existing studies on hyperbolic conservation laws with singular sources assume that the singular source depends continuously on the local state of the flow field, such as \cite{thanh2009riemann,greenberg1996well,das2020sharp,yang2013discontinuous,jin2001steady}.
Since the states of the flow fields on either side of the singular source are generally discontinuous, we focus on a more general model in which the source depends explicitly on the states on either side in the present work.
A typical example is the model of condensation of wet gas in a nozzle, where the rate of condensation is determined only by the state of upstream flow (see \cite{schnerr2005unsteadiness,2002On,cheng2010condensation}).

For simplicity, we assume that the value of $S$ is proportional to the flux upstream of the origin.
\begin{equation}\label{equation: value of source}
	S(U_-,U_+)=\begin{cases}
		S(F(U_-))=diag(k_1,k_2,k_3)F(U_-),\ \text{if}\ u_->0,u_+>0\\
		S(F(U_+))=diag(k_1,k_2,k_3)F(U_+),\ \text{if}\ u_-<0,u_+<0\\
		0,\ else\\
	\end{cases}
\end{equation}
where
\[ U_-=U(0-,t),\ U_+=U(0+,t), \]
$diag(k_1,k_2,k_3)$ is the diagonal matrix with constant diagonal element $k_1$, $k_2$ and $k_3$.
The models in \cite{2002On} can be considered as a special case of the above model.

A basic assumption is that the fluid cannot lose more mass, momentum and energy than the corresponding total amount flowing into the origin, which implies
\begin{equation*}
	k_1>-1,k_2>-1,k_3>-1.
\end{equation*}
When the fluid at the origin is in equilibrium, the upstream and downstream velocities have the same sign. For the non-equilibrium state where these two velocities have different signs, we set the value of source to zero. We assume that the fluid is polytropic ideal and thermally ideal, whose equation of state is given by
\begin{equation*}
	R T=p/\rho,\	p=(\gamma-1)\rho e,\ \gamma>1, 
\end{equation*}
where $\gamma$ is the ratio of specific heats, $e$ is the internal energy, $T$ is the temperature and $R$ is the universal gas constant. The Riemann problem of (\ref{equ: governing equation}) is a Cauchy problem with piecewise constant initial values
\begin{equation}\label{equ: Riemann problem 1}
	\begin{cases}
		\partial _tU+\partial _xF(U)=\delta(x)S \\
		U(x,0)=\begin{cases}
			U_L,\ x<0 \\
			U_R,\ x>0
		\end{cases}
	\end{cases}
\end{equation}

A typical feature of the Riemann solution of hyperbolic conservation law with singular sources is that it contains a new kind of wave, which appears as a stationary discontinuity. 
Without a restriction on this new discontinuity, the Riemann solution is not unique(see \cite{greenberg1997analysis,greenberg1996well}), thus we need to first set up an admissibility condition to select the physical solution. 
It is very difficult to establish an admissibility condition from the global solution, and very limited work has been done, such as a minimum principle proposed by J.M.Greenberg et al.\cite{greenberg1997analysis} by a prior estimation. 
A simpler and more practical approach is to limit the state on both sides of this new discontinuity, that is, to restrict the steady solution near the stationary discontinuity rather than the global Riemann solution, as in the Lax shock inequality(see \cite{lax1957hyperbolic,lax1973hyperbolic}). 
The admissibility condition established in this way has been successfully applied to the shallow water equations(see \cite{lefloch2011godunov,lefloch2007riemann}), the equations of flows in a nozzle(see \cite{kroner2008minimum,kroner2005numerical,lefloch2003riemann,thanh2009riemann}) and the equations of flow with heat addition(see \cite{anderson1990modern,cheng2010condensation,delale2007condensation,2002On,schnerr2005unsteadiness}). 
We find that those admissibility condition for different equations are equivalent in the representation of eigenvalues, which motivates us to propose an eigenvalue-based condition in this work. 
We apply the present proposed condition to (\ref{equ: governing equation}) and study the wave curve for the new disontinuity.
Following the nomenclature in \cite{kroner2008minimum,kroner2005numerical,lefloch2003riemann,lefloch2011godunov,lefloch2007riemann,thanh2009riemann,thanh2016properties}, the present proposed condition is still named as a monotonicity criterion, since it also leads to a kind of monotonicity of the new discontinuity.

In general, the singular source causes the structure of Riemann solution to become complicated. For example, the number of waves is not fixed, and the waves may interact with each other(see \cite{alcrudo2001exact,bernetti2008exact,cheng2010condensation,lefloch2011godunov,lefloch2007riemann,pares2019riemann,thanh2009riemann,yuTARiemann,yu2021waves}).
A common way to analyze such a kind of Riemann problem is to convert the original non-homogeneous system into a homogeneous system by adding an augmented equation, and the discontinuity induced by the singular sources can then be associated with a new linearly degenerate characteristic domain of augmented equations.
There are a number of representative works that analyse Riemann problem by the augmented equations, and we refer readers to \cite{bernetti2008exact,dumbser2013diffuse,goatin2004riemann,lefloch2011godunov,pares2019riemann,thanh2009riemann} for details.
In form, we can also convert the Riemann problem (\ref{equ: Riemann problem 1}) into the Riemann problem of a homogeneous system:
\begin{equation}\label{equ: augmentation}
	\begin{cases}
		\partial _t\begin{pmatrix}
			U\\h
		\end{pmatrix}
		+\partial _x\begin{pmatrix}
			F(U)\\0
		\end{pmatrix}
		-\begin{pmatrix}
			S\partial_xh\\
			0
		\end{pmatrix}=0, \\
		U(x,0)=\begin{cases}
			U_L,\ x<0 \\
			U_R,\ x>0
		\end{cases}\ 
		h(x,0)=\begin{cases}
			0,\ x<0\\
			1,\ x>0
		\end{cases}
	\end{cases}
\end{equation}
The augmented system is nonconservative equations, and the term $S\partial_xh$ is the product of two distributions at the origin. 
The theory of nonconservative product has been developed to give a clear definition of this type of product, see \cite{abgrall2010comment,castro2007well,gosse2000well,gosse2001well,greenberg1996well,le1989shock,pares2006numerical,1994Why}. 
However, we would not use (\ref{equ: augmentation}) to obtain the Riemann solution in this paper. 
The reason being that the discontinuity of $S$, which is different from other common singular sources, forces us to redesign reasonable non-conservative paths. 
In addition, Abgrall, Bacigaluppi and Tokareva (2018) commented that "A nonconservative formulation of the system and getting the correct solution has been a long-standing debate"(p.10)\cite{abgrall2018high}.

In this paper we obtain the exact solution to Riemann problem (\ref{equ: Riemann problem 1}) by a double classical Riemann problem(CRP) framework\cite{yuTARiemann,yu2021waves}. 
In the x-t plane, the singular source is constantly located at the t-axis, so that the Riemann solutions in the left and right half-planes of the origin satisfy the classical Euler equations. 
Under the double CRP framework, the solution of CRP in the left half-plane, the steady solution near the singular source and the solution of CRP in the right half-plane are discussed separately and then coupled to obtain the Riemann solution in the full plane.

This paper is organized as follows. 
In Section\ref{Elementary waves} we will present all the elementary waves of Riemann solution, including the three simple waves associated with the classical Euler equations and a new kind of discontinuity induced by singular sources. 
We will give the wave curve of that new discontinuity under the monotonicity criterion.
Section\ref{Structures of Riemann solution} contains the proof of all possible structures of  Riemann solution.
Finally, some conclusions will be given in Section\ref{Conclusions}.

\section{Elementary waves}\label{Elementary waves}
In this section we present elementary waves in self-similar solutions of Riemann problem (\ref{equ: Riemann problem 1}), including a new discontinuity induced by the singular source term, called the stationary wave, and elementary waves in the Riemann solutions of classical Euler equations.

\subsection{Stationary wave}

The singular source term causes the state of flow to jump when it passes through the origin. This discontinuity is called the stationary wave and is regarded as a new elementary wave of the Riemann solution. Its left-hand state $U_-$ and right-hand state $U_+$ satisfy the jump relation:
\begin{equation}\label{equ: jump relation 1}
	F(U_-)+S(F(U_-))=F(U_+),
\end{equation}
and the component-wise form is
\begin{equation}\label{equ: jump relation 2}
	\begin{cases}
		\rho_-u_-(1+k_1)=\rho_+u_+,\\
		(\rho_-u_-^2+p_-)(1+k_2)=\rho_+u_+^2+p_+,\\
		(E_-+p_-)u_-(1+k_3)=(E_++p_+)u_+.
	\end{cases}
\end{equation}
If the left-hand state $U_-$ and right-hand state $U_+$ of $U$ at the origin satisfy (\ref{equ: jump relation 1}), then $U$ is a steady solution of (\ref{equ: governing equation}) near the origin. $U_-$ and $U_+$ are called equilibrium states. By the first eqution of (\ref{equ: jump relation 2}), we know that $u_-$ and $u_+$ have the same sign. One of the values of $u_-$ and $u_+$ being zero results in the other also being zero, which causes the source term to vanish and the governing equation (\ref{equ: governing equation}) to degenerate to the classical Euler equations. Hence we will only consider the case of $u_->0,u_+>0$, as the case of $u_-<0,u_+<0$ can be changed to that by a transformation $x\mapsto -x,u\mapsto -u$.

When $S=0(k_1=k_2=k_3=0)$, there are two $U_+$ that satisfy (\ref{equ: jump relation 1}) for a given $U_-$, and the Mach numbers of $U_+$ are
\begin{equation*}
	M_+=M_-\ \text{or}\ 	\sqrt{\frac{(\gamma-1)M_-^2+2}{2\gamma M_-^2+1-\gamma}}.
\end{equation*}

When $S\ne 0$, we can obtain an explicit expression of $U_+$ by a tedious computation similar to that used in \cite{yuTARiemann,anderson1990modern,cheng2010condensation}. There are also two solutions of $U_+$, and their Mach numbers are denoted by
\begin{equation}\label{equ:  two branches}
	M_+^{(1)}=\sqrt{\frac{1-I}{1+\gamma I}},\quad 
	M_+^{(2)}=\sqrt{\frac{1+I}{1-\gamma I}},
\end{equation}
where
\begin{equation*}
	I=I(M_-)=\frac{\sqrt{(\gamma M_-^2+1)^2-(\gamma+1)M_-^2[(\gamma-1)M_-^2+2](1+k_1)(1+k_3)/(1+k_2)^2}}{\gamma M_-^2+1}.
\end{equation*}
$M_+^{(1)}$ is well-defined iff
\begin{equation}\label{equ: Mach number condition1}
	(\gamma M_-^2+1)^2\geq (\gamma+1)M_-^2[(\gamma-1)M_-^2+2](1+k_1)(1+k_3)/(1+k_2)^2.
\end{equation}
$M_+^{(2)}$ is well-defined iff
\begin{equation}\label{equ: Mach number condition2}
	\text{equation}(\ref{equ: Mach number condition1})\ \&\ \gamma I\leq 1.
\end{equation}

Since $I$ is not less than zero, we have
\begin{equation}\label{equ:two Mach numbers and 1}
	M_+^{(1)}\leq 1,\ M_+^{(2)}\geq 1.
\end{equation}
$M_+^{(1)}$ and $M_+^{(2)}$ are referred to as the subsonic and supersonic branches, respectively. Besides, a standard computation(see \cite{anderson1990modern}) gives rise to
\begin{equation*}
	\begin{aligned}
		&\frac{\rho_+}{\rho_-}=\left(\frac{M_-}{M_+} \right) ^2\frac{\gamma M_+^2+1}{\gamma M_-^2+1}\frac{(1+k_1)^2}{1+k_2},\\
		&\frac{u_+}{u_-}=\left(\frac{M_+}{M_-} \right) ^2\frac{\gamma M_-^2+1}{\gamma M_+^2+1}\frac{1+k_2}{1+k_1},\\
		&\frac{p_+}{p_-}=\frac{\gamma M_-^2+1}{\gamma M_+^2+1}(1+k_2).
	\end{aligned}
\end{equation*}
That is, the ratio of other states between downstream and upstream can be determined by the Mach numbers $M_-$ and $M_+$. With $M_-\in (0,\infty)$ as an independent variable, $M_+^{(1)}$ and $M_+^{(2)}$ defined by (\ref{equ:  two branches}) can be regarded as functions of $M_-$, then it can be proved that
\begin{Property}\label{property: extreme value}
	Under (\ref{equ: Mach number condition1}), $M_+^{(1)}(M_-)$ is increasing when $M_-\in (0,1]$ and decreasing when $M_-\in[1,\infty)$. Under (\ref{equ: Mach number condition2}), $M_+^{(2)}(M_-)$ is decreasing when $M_-\in (0,1]$ and increasing when $M_-\in[1,\infty)$.
\end{Property}

The quasilinear form of the inhomogeneous Euler equations (\ref{equ: governing equation}) is
\begin{equation*}
	\partial _tU+A(U)\partial _xU=\delta(x)S,
\end{equation*}
where $A(U)=\partial F/ \partial U$. This equation is strictly hyperbolic, and the eigenvalues of $A(U)$ are
\[
\lambda_1=u-a<\lambda_2=u<\lambda_3=u+a,
\]
where $a$ is the speed of sound.

The two solutions (\ref{equ:  two branches}) of $M_+$ may be nonphysical, therefore an admissibility condition called monotonicity criterion is given below to select the physical solution of (\ref{equ: jump relation 1}). 
\begin{Monotonicity Criterion}\label{criterion1}
	Given coefficients $k_1$, $k_2$ and $k_3$, the left-hand state $U_-$ and right-hand state $U_+$ of a statiaonry wave satisfies $\lambda_k(U_-)\cdot \lambda_k(U_+)\geq 0,\forall 1\leq k\leq 3$.
\end{Monotonicity Criterion}

Criterion\ref{criterion1} restricts each corresponding eigenvalue of the left-hand and right-hand states of stationary wave from having opposite signs. When a certain eigenvalue of the state on one side is zero, Criterion\ref{criterion1} loses its restriction on the corresponding eigenvalue of the state on the other side. Following the definition of thermal choke in \cite{2002On,delale2007condensation,cheng2010condensation,anderson1990modern}, we define this particular state as choke as well.
\begin{Definition}
	A stationary wave is said to be choked if its left-hand state $U_-$ and right-hand state $U_+$ satisfy $\exists 1\leq k\leq m, s.t.\ \lambda_k(U_-)\cdot \lambda_k(U_+)=0$.
\end{Definition}
A Riemann solution is said to be choked if it contains a choked stationary wave. 

\begin{Remark}
	Criterion\ref{criterion1} can be applied to Generalized hyperbolic conservation laws with singular source terms:
	\begin{equation}\label{equ: generalized equation}
		\partial _t U+\partial _x F(U)=\delta (x) S,
	\end{equation}
	where $U$, $F$ and $S$ are all vector-valued functions of m components. If the self-similar solution to the Riemann problem of (\ref{equ: generalized equation}) contains a stationary wave on the t-axis with the left-hand state $U_-$ and the right-hand state $U_+$, then it satisfies $\lambda_k(U_-)\cdot \lambda_k(U_+)\geq 0,\forall 1\leq k\leq m$, where $\lambda_k(1\leq k\leq m)$ is the kth eigenvalue of the $m\times m$ matrix $\partial F/ \partial U$.  
\end{Remark}

\begin{Remark}
	Greenberg, Leroux, Baraille and Noussair\cite{greenberg1997analysis} studied the Riemann problem of the scalar equation with a singular source term
	\begin{equation}\label{equ: scalar equation}
		\begin{cases}
			\partial _t U+\partial _x F(U)=c \delta(x),\\
			U(x,0)=\begin{cases}
				U_L,\ x<0 \\
				U_R,\ x>0
			\end{cases}
		\end{cases}
	\end{equation} 
	where $U(x,t),F(U)$ are both scalar-valued functions, $F(0)=0,F''>0$, and $c$ is a constant coefficient. Greenberg et al. proposed a minimum principle, which is based on a priori estimate of the integral equation, to select the physical Riemann solutions. The expression of Criterion\ref{criterion1} applied to (\ref{equ: scalar equation}) is $U_-\cdot U_+\geq 0$, by which we can obtain the values of $U_-$ and $U_+$ of Riemann solution. It's easy to verify that the Riemann solutions based on Greenberg's principle and the Riemann solutions based on the present Monotonicity Criterion are exactly the same in all cases.
\end{Remark}

\begin{Remark}
	The shallow water equations with discontinuous topography is a widely studied hyperbolic conservation law with singular source terms, and its Riemann problem is
	\begin{equation*}
		\begin{cases}
			\frac{\partial h}{\partial t}+\frac{\partial hu}{\partial x}=0,\\
			\frac{\partial hu}{\partial t}+\frac{\partial (h(u^2+gh/2))}{\partial x}=-gh\frac{\partial a}{\partial x},\\
			(h,u,a)(x,0)=\begin{cases}
				(h_L,u_L,a_L),\quad x<0\\
				(h_R,u_R,a_R),\quad x>0
			\end{cases}
		\end{cases}
	\end{equation*}	
	Here, $h$ denotes the height of the water from the bottom to the surface, $u$ the velocity of the fluid, $g$ the gravity constant, and $a$ the height of the river bottom from a given level. The source term of shallow water equations differs from (\ref{equ: governing equation}) in that it is continuously dependent on the conservative variable, therefore additional approximations are required to give a jump relation. The jump relation used in \cite{lefloch2007riemann,lefloch2011godunov,alcrudo2001exact} is
	\begin{equation*}
		\begin{cases}
			h_-u_-=h_+u_+,\\
			h_-+\frac{u_-^2}{2g}+a_L=h_++\frac{u_+^2}{2g}+a_R,
		\end{cases}
	\end{equation*}
	The subscripts "\_" and "+" still denote the left-hand and right-hand states of t-axis of Riemann solution. The expression for Criterion\ref{criterion1} applied to the shallow water equations is $F_-\cdot F_+\geq 0$, where $F$ is Froude number, which can be proved to be equivalent to the admissible principles used in previous literature(see \cite{lefloch2007riemann,lefloch2011godunov,thanh2016properties}).
\end{Remark}

\begin{Remark}
	The jump relation (\ref{equ: jump relation 1}) with $k_1=k_2=0,k_3>0$ is called equilibrium equation of heat addition problem(see \cite{schnerr2005unsteadiness,2002On,delale2007condensation}). The Euler equations with this particular source term has been used to study the condensation-induced waves in a slender Laval nozzle in \cite{cheng2010condensation}. The physical properties force the change of the state during heating not to cross the maximum entropy point in the Rayleigh line (see \cite{anderson1990modern}), which equates to the Criterion\ref{criterion1}.
\end{Remark}

We call $U_+$ the solution to a stationary wave if $U_+$ satisfies (\ref{equ: jump relation 1}) and Criterion\ref{criterion1} for given $U_-$ and coefficients $k_1,k_2,k_3$. According to Criterion\ref{criterion1} and (\ref{equ:two Mach numbers and 1}), we can select the physical solutions between the two branches $M_+^{(1)}$ and $M_+^{(2)}$:
\begin{equation}\label{equ: downstream Mach number}
	M_+=\begin{cases}
		M_+^{(1)}=\sqrt{\frac{1-I}{1+\gamma I}},\ \text{if}\ M_-\leq 1\\
		M_+^{(2)}=\sqrt{\frac{1+I}{1-\gamma I}},\ \text{if}\ M_-\geq 1
	\end{cases}
\end{equation}
According to (\ref{equ: downstream Mach number}), (\ref{equ: Mach number condition1}) and (\ref{equ: Mach number condition2}), the range of $M_-$ that ensures the existence of the physical solutions is
\begin{equation*}
	\begin{aligned}
		&\Gamma_{-,k_1,k_2,k_3}\mathop{=}\limits^{def}\left( (0,1] \cup \left( [1,+\infty) \cap \left\lbrace M_-|\gamma I(M_-)\leq 1 \right\rbrace  \right)  \right) \cap\\
		&\left\lbrace M_-|(\gamma M_-^2+1)^2\geq (\gamma+1)M_-^2[(\gamma-1)M_-^2+2](1+k_1)(1+k_3)/(1+k_2)^2 \right\rbrace ,
	\end{aligned}
\end{equation*}
which will be written simply as $\Gamma_-$ when no confusion can arise. We can now give a clear definition of the stationary wave.
\begin{Proposition}[Stationary wave curves]\label{proposition: wave curve}
	Given a left-hand state $U_-$ with $u_->0$, the stationary wave curve $\mathcal{D}(U_-,k_1,k_2,k_3)$ consisting of all right-hand states $U_+$ that can be connected to $U_-$ by a stationary wave with given coefficients $k_1,k_2,k_3$ under Criterion\ref{criterion1}, is
	\begin{equation*}
		\begin{aligned}
			\mathcal{D}(U_-,k_1,k_2,k_3):\quad
			&M_+=\begin{cases}
				M_+^{(1)}=\sqrt{\frac{1-I}{1+\gamma I}},\ M_-\in \Gamma_-\cap (0,1]\\
				M_+^{(2)}=\sqrt{\frac{1+I}{1-\gamma I}},\ M_-\in \Gamma_-\cap [1,+\infty)
			\end{cases}\\
			&\rho_+=\rho_-\left(\frac{M_-}{M_+} \right) ^2\frac{\gamma M_+^2+1}{\gamma M_-^2+1}\frac{(1+k_1)^2}{1+k_2},\\
			&u_+=u_-\left(\frac{M_+}{M_-} \right) ^2\frac{\gamma M_-^2+1}{\gamma M_+^2+1}\frac{1+k_2}{1+k_1},\\
			&p_+=p_-\frac{\gamma M_-^2+1}{\gamma M_+^2+1}(1+k_2).
		\end{aligned}
	\end{equation*}
\end{Proposition}

\begin{Note}\label{note: stationary wave}
	Although we refer to the set $\mathcal{D}(U_-,k_1,k_2,k_3)$ of solutions as a wave curve, it is not a one-parameter family. If the coefficients $k_1,k_2,k_3$ are taken as parameters, it is actually a three-parameter family. For given parameters $k_1,k_2,k_3$, each left-hand state $U_-$ determines its corresponding right state $U_+$, and we therefore say that the stationary wave has no degrees of freedom.
\end{Note}

\begin{Property}\label{property: monotonicity}
	The $M_+$ defined by (\ref{equ: downstream Mach number}) is a monotone increasing function of $M_-\in \Gamma_-\cap (0,1]$ and $M_-\in \Gamma_-\cap [1,+\infty)$.
\end{Property}
It is clear from PropertyProperty\ref{property: extreme value} that both branches $M_+^{(1)},M_+^{(2)}$ as a function of $M_-$ have a unique extreme value at $M_-=1$, and PropertyProperty\ref{property: monotonicity} shows that the Criterion\ref{criterion1} works by choosing a branch such that the Mach number $M_+$ of physical solution is a piecewise monotonic function of $M_-$.

The solution of a stationary wave depends on the solution of Mach number $M_+$, as each $M_+$ corresponds to a unique state $U_+$. Similarly, we denote by $\Gamma_{+,k_1,k_2,k_3}$ the set of $M_+$, which is also written as $\Gamma_+$ when no confusion can arise. As will be seen in the next section, $\Gamma_-$ and $\Gamma+$ significantly affect the structure of the Riemann solution, therefore we give more explicit expressions for $\Gamma_-$ and $\Gamma_+$ below.

With the specific eigenvalues of Euler equations, an equivalent expression for Criterion\ref{criterion1} as applied to one-dimensional Euler equations is
\begin{Monotonicity Criterion}\label{criterion2}
	Given coefficients $k_1$, $k_2$ and $k_3$, the left-hand state $U_-$ and right-hand state $U_+$ of a stationary wave satisfies $u_-\cdot u_+\geq 0$ and $(M_--1)\cdot(M_+-1)\geq 0$.
\end{Monotonicity Criterion}

\begin{figure}[h]
	\label{figure: admissible region}
	\centering
	\includegraphics[width=0.4\linewidth]{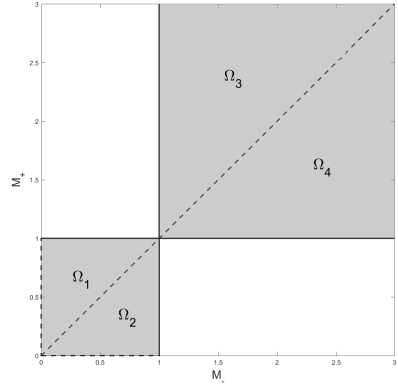}
	\caption{The admissible region of monotonicity criterion in $M_--M_+$ graph.}
\end{figure}

The admissible region of Criterion\ref{criterion2} in $M_--M_+$ plane is $(0,1]\times (0,1] \cup [1,+\infty)\times [1,+\infty) $, as shown in the gray area in Figure\ref{figure: admissible region}.
We divide this area into four sub-areas as follows.
\begin{equation*}
	\begin{aligned}
		&\Omega_1=\left\lbrace (M_-,M_+)| 0<M_-\leq 1, 0<M_+\leq 1, M_+>M_- \right\rbrace \\
		&\Omega_2=\left\lbrace (M_-,M_+)| 0<M_-\leq 1, 0<M_+\leq 1, M_+<M_- \right\rbrace \\
		&\Omega_3=\left\lbrace (M_-,M_+)| M_-\geq 1, M_+\geq 1, M_+>M_- \right\rbrace \\
		&\Omega_4=\left\lbrace (M_-,M_+)| M_-\geq 1, M_+\geq 1, M_+<M_- \right\rbrace 
	\end{aligned}
\end{equation*}

\begin{Property}\label{property: three cases}
	For given $U_-$ and coefficients $k_1,k_2,k_3$, the $U_+\in D(U_-,k_1,k_2,k_3)$ satisfies
	\begin{itemize}
		\item[(i)] if $k>0$, then $|M_+-1|<|M_--1|$;
		\item[(ii)] if $-1<k<0$, then $|M_+-1|>|M_--1|$;
		\item[(iii)] if $k=0$, then $|M_+-1|=|M_--1|$,
	\end{itemize}
	where $k=(1+k_3)(1+k_1)/(1+k_2)^2-1$.
\end{Property}
The proof is straightforward. Property\ref{property: three cases} shows that the physical solution of $M_+$ is located in $\Omega_1 \cup \Omega_4$ when $k>0$, and is located in $\Omega_2 \cup \Omega_3$ when $k<0$.

\begin{figure}[htbp]\label{figure: Mach numbers}
	\centering
	\subfigure{
		\includegraphics[width=0.4\linewidth]{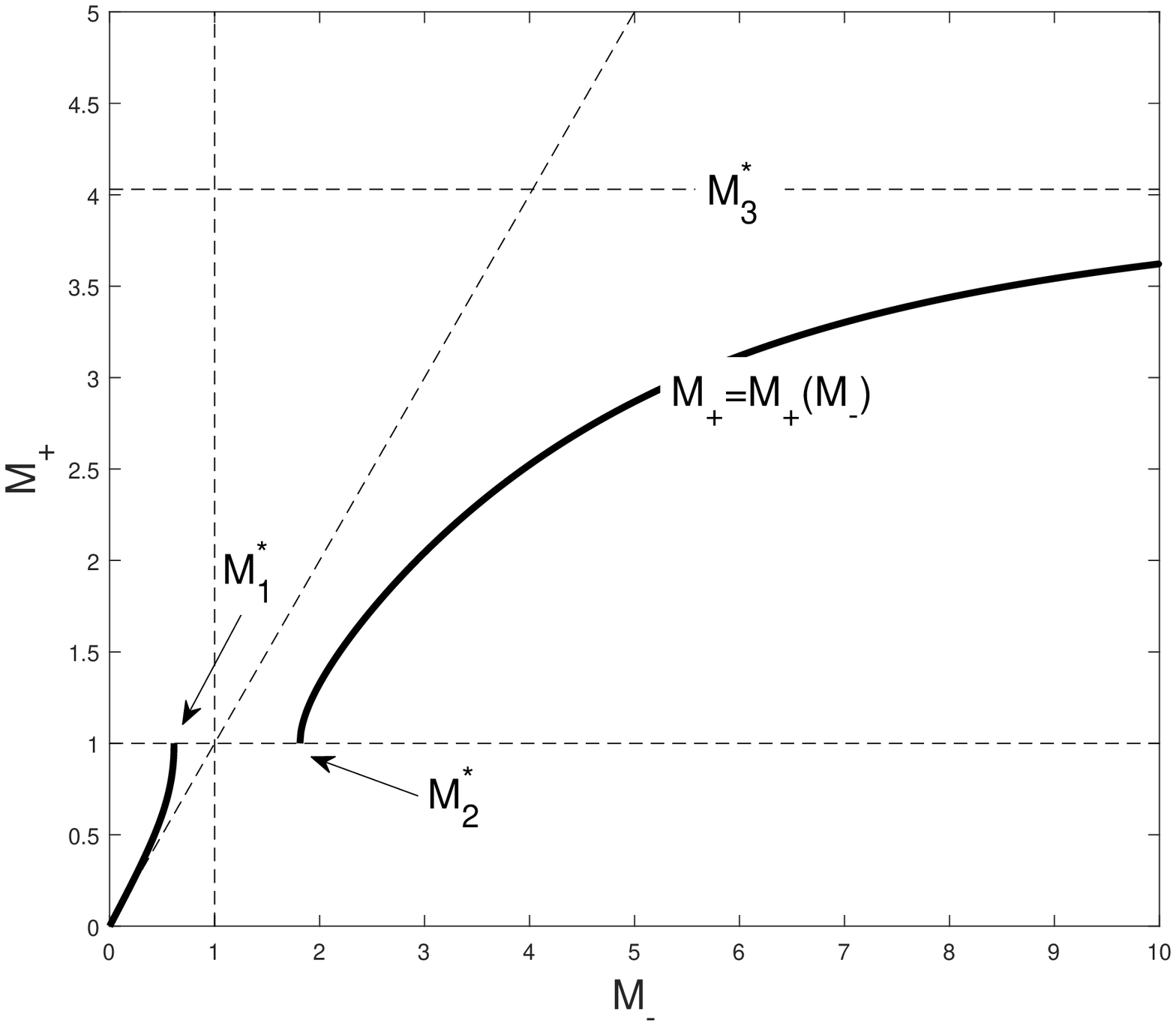}
	}
	\quad
	\subfigure{
		\includegraphics[width=0.42\linewidth]{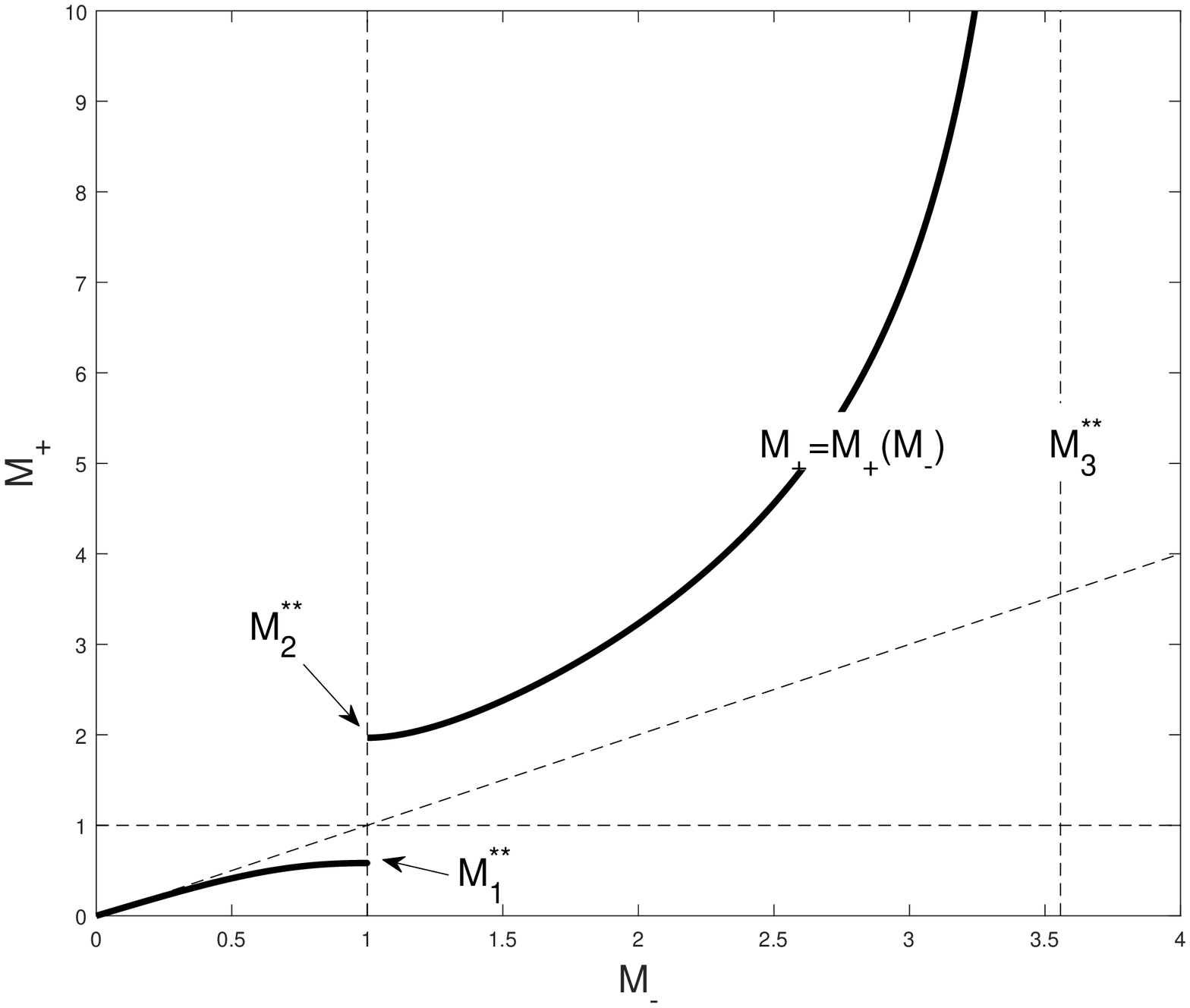}
	}
	\caption{The function $M_+=M_+(M_-)$ defined by (\ref{equ: downstream Mach number}). Left: $\gamma=1.4,k=0.2$. Right: $\gamma=1.4,k=-0.2$.}
\end{figure}

We first consider the case of $k>0$. A typical example of the function $M_+(M_-)$ for $k>0$ is shown in the left diagram in Figure\ref{figure: Mach numbers}. When $M_-\in \Gamma_-\cap (0,1]$, simple computation gives rise to
\begin{equation*}
	\begin{aligned}
		&\lim_{M_-\to 0}M_+(M_-)=0,\\
		&\left( M_+^{(1)}\right) ^{-1}(1)=\begin{cases}
			\sqrt{\frac{k\gamma+k+1- (\gamma+1)\sqrt{k(k+1)}}{k+1-k\gamma^2}}<1,\ &\text{if}\ k>0,\ k \ne \frac{1}{\gamma^2-1}\\
			\sqrt{\frac{\gamma-1}{2\gamma}}<1,\ &\text{if}\ k = \frac{1}{\gamma^2-1}
		\end{cases}. 
	\end{aligned}
\end{equation*}
When $M_-\in \Gamma_-\cap [1,+\infty)$, we assume that $k<\frac{1}{\gamma^2-1}$, since $M_+$ is unsolvable when $k\geq \frac{1}{\gamma^2-1}$. It follows that
\begin{equation*}
	\begin{aligned}
		&\lim_{M_-\to +\infty}M_+(M_-)=\sqrt{\frac{\gamma+\sqrt{1-k(\gamma^2-1)}}{\gamma-\gamma\sqrt{1-k(\gamma^2-1)}}}>1\\
		&\left( M_+^{(2)}\right) ^{-1}(1)=\sqrt{\frac{k\gamma+k+1+ (\gamma+1)\sqrt{k(k+1)}}{k+1-k\gamma^2}}> 1
	\end{aligned}
\end{equation*}

We use the following notations to denote these critical Mach numbers:
\begin{equation}\label{equ: critical Mach number definition1}
	\begin{aligned}
		&M_1^*=\begin{cases}
			\sqrt{\frac{k\gamma+k+1- (\gamma+1)\sqrt{k(k+1)}}{k+1-k\gamma^2}},\ &\text{if}\ k \ne \frac{1}{\gamma^2-1},\ k>0\\
			\sqrt{\frac{\gamma-1}{2\gamma}},\ &\text{if}\ k = \frac{1}{\gamma^2-1}
		\end{cases} \\
		&M_2^*=\begin{cases}
			\sqrt{\frac{k\gamma+k+1+ (\gamma+1)\sqrt{k(k+1)}}{k+1-k\gamma^2}},\ &\text{if}\ 0<k< \frac{1}{\gamma^2-1}\\
			+\infty, \ &\text{if}\ k\geq \frac{1}{\gamma^2-1}
		\end{cases}\\
		&M_3^*=\begin{cases}
			\sqrt{\frac{\gamma+\sqrt{1-k(\gamma^2-1)}}{\gamma-\gamma\sqrt{1-k(\gamma^2-1)}}}, \ &\text{if}\ 0<k< \frac{1}{\gamma^2-1}\\
			+\infty, \ &\text{if}\ k\geq \frac{1}{\gamma^2-1}
		\end{cases}
	\end{aligned}
\end{equation}
If $k\ne 0$, then $M_1^*< M_2^*$, which implies that the solution of a stationary wave may not exist. To be precise, the left-hand state $U_-$ and right-hand state $U_+$ of a stationary wave should satisfy
\begin{equation*}
	M_-\in (0,M_1^*]\&M_+\in (0,1] \ \text{or } M_-\in [M_2^*,\infty)\&M_+\in [1,M_3^*). 
\end{equation*}

Secondly, we consider the case of $-1<k<0$. A typical example of the function $M_+(M_-)$ for $-1<k<0$ is shown in the right diagram in Figure\ref{figure: Mach numbers}. It is obvious that (\ref{equ: jump relation 2}) is equivalent to
\begin{equation}\label{equ: transformation}
	\begin{aligned}
		&\begin{cases}
			\rho_-u_-=\rho_+u_+\left( 1+\frac{-k_1}{1+k_1} \right), \\
			\rho_-u_-^2+p_-=(\rho_+u_+^2+p_+)\left( 1+\frac{-k_2}{1+k_2} \right),\\
			(E_-+p_-)u_-=(E_++p_+)u_+\left( 1+\frac{-k_3}{1+k_3} \right),
		\end{cases}\\
		&\frac{\left( 1+\frac{-k_1}{1+k_1}\right) \left( 1+\frac{-k_3}{1+k_3} \right)} { \left( 1+\frac{-k_2}{1+k_2} \right)^2 }-1
		=\frac{-k}{1+k}.
	\end{aligned}
\end{equation}
Criterion\ref{criterion2} is symmetric for $U_-$ and $U_+$, so we are able to use all the conclusions of case $k>0$ directly with the help of transformation (\ref{equ: transformation}). Set
\begin{equation}\label{equ: critical Mach number definition2}
	\begin{aligned}
		&M_1^{**}=\sqrt{\frac{1-\sqrt{-k}}{1+\gamma\sqrt{-k}}},\\
		&M_2^{**}=\begin{cases}
			\sqrt{\frac{1+\sqrt{-k}}{1-\gamma\sqrt{-k}}},\ &\text{if}\  -\frac{1}{\gamma^2}<k<0\\
			+\infty, \ &\text{if}\ -1<k\leq -\frac{1}{\gamma^2}
		\end{cases}\\
		&M_3^{**}=\begin{cases}
			\sqrt{\frac{\gamma\sqrt{1+k}+\sqrt{1+k\gamma^2}}{\gamma\sqrt{1+k}-\gamma\sqrt{1+k\gamma^2}}}, \ &\text{if}\ -\frac{1}{\gamma^2}<k<0\\
			1, \ &\text{if}\ -1<k\leq -\frac{1}{\gamma^2}
		\end{cases}
	\end{aligned}
\end{equation}
A trivial verification shows that when $-1<k<0$ the left-hand state $U_-$ and right-hand state $U_+$ should satisfy
\begin{equation*}
	M_-\in (0,1]\&M_+\in (0,M_1^{**}] \ \text{or } M_-\in [1,M_3^{**})\&M_+\in [M_2^{**},+\infty). 
\end{equation*}

Finally,  the range of left-hand Mach number and right-hand Mach number for the case $k=0$ is
\begin{equation*}
	M_+=M_-\in(0,+\infty).
\end{equation*}

For given $k_1,k_2,k_3>-1$, the sets of left-hand Mach number $M_-$ and right-hand Mach number $M_+$ are shown in table\ref{table: left and right Mach numbers}.
In addition to $M_-=1$, we can choose the solution branch of stationary wave according to $M_-$. When $M_-=1$, it can be shown that $M_+^{(1)}(1)=M_+^{(2)}(1)$ when and only when $k=0$, which means that the solution to stationary wave is unique when $k=0$, but not when $k< 0$.

\begin{Property}\label{property: existence and uniqueness}
	For the stationary wave with given left-hand state $U_-$ and coefficients $k_1,k_2,k_3$,
	\begin{itemize}
		\item[(i)] the solution does not exist iff
		\begin{equation*}
			k>0\ \&\ M_1^{*}<M_-<M_2^{*}\quad \text{or}\quad 
			k<0\ \&\ M_-\geq M_3^{**};
		\end{equation*}
		\item [(ii)] the solution is double iff $k<0\ \&\ M_-=1$;
		\item[(iii)] the solution exists and is unique in other cases.
	\end{itemize}
\end{Property}

Here is a practical definition of stationary wave.
\begin{Proposition}[Stationary wave curves]\label{proposition: stationary wave curve1}
	Given a left-hand state $U_-$, the stationary wave curve $\mathcal{D}(U_-,k_1,k_2,k_3)$ consisting of all right-hand states $U_+$ that can be connected to $U_-$ by the stationary wave with coefficients $k_1,k_2,k_3$ which satisfies Criterion\ref{criterion1}, is
	\begin{equation*}
		\begin{aligned}
			D(U_-,k_1,k_2,k_3):\quad
			&\text{if}\ k>0,\ M_+=\begin{cases}
				\sqrt{\frac{1-I}{1+\gamma I}},\ 0<M_-\leq M_1^{*}\\
				\sqrt{\frac{1+I}{1-\gamma I}},\ M_2^{*}\leq M_-<+\infty
			\end{cases}\\
			&\text{if}\ k<0,\ M_+=\begin{cases}
				\sqrt{\frac{1-I}{1+\gamma I}},\ 0<M_-\leq 1\\
				\sqrt{\frac{1+I}{1-\gamma I}},\ 1\leq M_-< M_3^{**}\\
			\end{cases}\\
			&\text{if}\ k=0,\ M_+=\begin{cases}
				\sqrt{\frac{1-I}{1+\gamma I}},\ 0<M_-\leq 1\\
				\sqrt{\frac{1+I}{1-\gamma I}},\ 1<M_-<+\infty\\
			\end{cases}\\			
			&\rho_+=\rho_-\left(\frac{M_-}{M_+} \right) ^2\frac{\gamma M_+^2+1}{\gamma M_-^2+1}\frac{(1+k_1)^2}{1+k_2},\\
			&u_+=u_-\left(\frac{M_+}{M_-} \right) ^2\frac{\gamma M_-^2+1}{\gamma M_+^2+1}\frac{1+k_2}{1+k_1},\\
			&p_+=p_-\frac{\gamma M_-^2+1}{\gamma M_+^2+1}(1+k_2).
		\end{aligned}
	\end{equation*}
	Here, 
	\[
	I=\frac{\sqrt{(\gamma M_-^2+1)^2-(\gamma+1)M_-^2[(\gamma-1)M_-^2+2](1+k)}}{\gamma M_-^2+1},
	\]
	$k=(1+k_1)(1+k_3)/(1+k_2)^2-1$ and $M_1^{*},M_2^{*},M_3^{**}$ are defined by (\ref{equ: critical Mach number definition1}) and (\ref{equ: critical Mach number definition2}).
\end{Proposition}

\begin{table}[htbp]
	\centering
	\caption{The $\Gamma_{-,k_1,k_2,k_3}$ and $\Gamma_{+,k_1,k_2,k_3}$ for different $k$.}
	\label{table: left and right Mach numbers}
	\begin{tabular}{ccccc}
		\toprule
		\specialrule{0em}{3pt}{2pt}
		\multirow{2}*{$k$}&\multicolumn{2}{c}{subsonic branch}&\multicolumn{2}{c}{supersonic branch}\\
		\specialrule{0em}{1pt}{2pt}
		&$\Gamma_{-}$& $\Gamma_{+}$&$\Gamma_{-}$& $\Gamma_{+}$\\
		\specialrule{0em}{3pt}{2pt}	
		\midrule
		\specialrule{0em}{1pt}{2pt}
		$(0,+\infty)$&$(0,M_1^*]$&$(0,1]$&$\left[M^*_2,+\infty \right) $&$\left[1,M_3^{*} \right) $\\
		\specialrule{0em}{1pt}{2pt}
		${0}$&$(0,1]$&$(0,1]$&$\left[1,+\infty \right) $&$\left[1,+\infty \right) $\\
		\specialrule{0em}{1pt}{2pt}
		$(-\infty,0)$&$(0,1]$&$(0,M_1^{**}]$&$\left[1,M_3^{**} \right)$&$\left[M_2^{**},+\infty \right)$\\
		\specialrule{0em}{3pt}{2pt}
		\bottomrule
	\end{tabular}
\end{table}

By the transformation (\ref{equ: transformation}), we can directly obtain the backward wave curve of a stationary wave.
\begin{Proposition}[Stationary wave curves]\label{proposition: stationary wave curve2}
	Given a right-hand state $U_+$, the stationary wave curve $\mathcal{D}^B(U_+,k_1,k_2,k_3)$ consisting of all left-hand states $U_-$ that can be connected to $U_+$ by the stationary wave with coefficients $k_1,k_2,k_3$ which satisfies Criterion\ref{criterion1}, is
	\begin{equation*}
		\begin{aligned}
			D^B(U_+,k_1,k_2,k_3):\quad
			&\text{if}\ k>0,\ M_-=\begin{cases}
				\sqrt{\frac{1-I^B}{1+\gamma I^B}},\ 0<M_+\leq 1\\
				\sqrt{\frac{1+I^B}{1-\gamma I^B}},\ 1\leq M_+< M^{*}_3\\
			\end{cases}\\
			&\text{if}\ k<0,\ M_-=\begin{cases}
				\sqrt{\frac{1-I^B}{1+\gamma I^B}},\ 0<M_+\leq M^{**}_1\\
				\sqrt{\frac{1+I^B}{1-\gamma I^B}},\ M^{**}_2\leq M_+<+\infty
			\end{cases}\\
			&\text{if}\ k=0,\ M_-=\begin{cases}
				\sqrt{\frac{1-I^B}{1+\gamma I^B}},\ 0<M_-\leq 1\\
				\sqrt{\frac{1+I^B}{1-\gamma I^B}},\ 1<M_-<+\infty\\
			\end{cases}\\			
			&\rho_-=\rho_+\left(\frac{M_+}{M_-} \right) ^2\frac{\gamma M_-^2+1}{\gamma M_+^2+1}\frac{1+k_2}{(1+k_1)^2},\\
			&u_-=u_+\left(\frac{M_-}{M_+} \right) ^2\frac{\gamma M_+^2+1}{\gamma M_-^2+1}\frac{1+k_1}{1+k_2},\\
			&p_-=p_+\frac{\gamma M_+^2+1}{\gamma M_-^2+1}\frac{1}{1+k_2}.
		\end{aligned}
	\end{equation*}
	Here, 
	\[
	I^B=\frac{\sqrt{(\gamma M_+^2+1)^2-(\gamma+1)M_+^2[(\gamma-1)M_+^2+2]/(1+k)}}{\gamma M_+^2+1},
	\]
	$k=(1+k_1)(1+k_3)/(1+k_2)^2-1$ and $M_1^{**},M_2^{**},M_3^{*}$ are defined by (\ref{equ: critical Mach number definition1}) and (\ref{equ: critical Mach number definition2}). 
\end{Proposition}

By the above discussion we see that
\begin{Property}
	The $M_+$ defined by (\ref{equ: downstream Mach number}) is a global monotone increasing function of $M_-$.
\end{Property}

\subsection{Shock waves, rarefaction waves and contact discontinuities}

In this section we will introduce the elementary waves other than stationary waves, including Lax shock waves, rarefaction waves and contact discontinuities. These three elementary waves are associated with the characteristic domains of matrix $A=\partial F/\partial U$ and therefore referred to as classical elementary waves. Here we only list some conclusions to be used in this paper. For detailed introductions, we refer the reader to \cite{toro2013riemann,chang1989riemann}. 

The characteristic domains of $\lambda_1=u-a$ and $\lambda_3=u+a$ are both genuinely nonlinear. The elementary wave associated with these two characteristic domains is either a shock wave or a rarefaction wave. A Lax shock wave is the discontinuity travelling at the speed $\sigma$, the left-hand and right-hand states of which are denoted as $U_1$ and $U_2$, respectively. 
Subject to the Lax shock inequality(\cite{lax1957hyperbolic,lax1973hyperbolic})
\begin{equation*}
	\begin{aligned}
		&\text{1-shock: }\lambda_1(U_1)>\sigma>\lambda_1(U_2),\\
		&\text{3-shock: }\lambda_3(U_1)>\sigma>\lambda_3(U_2),
	\end{aligned}
\end{equation*}
we can obtain the wave curve by the Rankine-Hugoniot condition:
\begin{equation*} 
	F(U_2)-F(U_1)=\sigma (U_2-U_1).
\end{equation*} 

\begin{Proposition}[Shock wave curves]
	Given a left-hand state $U_0$, the 1-shock curve $\mathcal{S}_1(U_0)$ consisting of all right-hand states $U$ that can be connected to $U_0$ by a Lax shock is
	\begin{equation*}
		\begin{aligned}
			\mathcal{S}_1(U_0):\ &\rho=\rho(p,U_0)=\rho_0\frac{(\gamma-1)p_0+(\gamma+1)p}{(\gamma-1)p+(\gamma+1)p_0},\ p\geq p_0\\
			&u=u(p,U_0)=u_0-(p-p_0)\sqrt{\frac{1}{\rho_0}\frac{2}{(\gamma+1)p+(\gamma-1)p_0}},\ p\geq p_0
		\end{aligned}
	\end{equation*}
	Given a right-hand state $U_0$, the 3-shock curve $\mathcal{S}_3(U_0)$ consisting of all left-hand states $U$ that can be connected to $U_0$ by a Lax shock is
	\begin{equation*}
		\begin{aligned}
			\mathcal{S}_3(U_0):\ &\rho=\rho(p,U_0)=\rho_0\frac{(\gamma-1)p_0+(\gamma+1)p}{(\gamma-1)p+(\gamma+1)p_0},\ p\geq p_0\\
			&u=u(p,U_0)=u_0+(p-p_0)\sqrt{\frac{1}{\rho_0}\frac{2}{(\gamma+1)p+(\gamma-1)p_0}},\ p\geq p_0
		\end{aligned}
	\end{equation*}
\end{Proposition}

A rarefaction wave is a piecewise smooth self-similar solution. The Riemann invariants across the 1-rarefaction and 3-rarefaction are 
\begin{align*}
	u+\frac{2a}{\gamma-1}=\text{constant},\ s=\text{constant}(\text{across}\ \lambda_1),\\
	u-\frac{2a}{\gamma-1}=\text{constant},\ s=\text{constant}(\text{across}\ \lambda_3),\\
\end{align*}
where $s$ is the entropy. Then we can get the wave curve of rarefaction waves.
\begin{Proposition}[Rarefaction wave curves]
	Given a left-hand state $U_0$, the 1-rarefaction curve $\mathcal{R}_1(U_0)$ consisting of all right-hand states $U$ that can be connected to $U_0$ by a rarefaction wave is
	\begin{equation*}
		\begin{aligned}
			\mathcal{R}_1(U_0):\ &\rho=\rho(p,U_0)=\rho_0\left( \frac{p}{p_0}\right) ^{\frac{1}{\gamma}},\ p\leq p_0\\
			&u=u(p,U_0)=u_0-\frac{2a_0}{\gamma-1}\left[ \left( \frac{p}{p_0}\right) ^{\frac{\gamma-1}{2\gamma}}-1\right] ,\ p\leq p_0\\
		\end{aligned}
	\end{equation*}
	Given a right-hand state $U_0$, the 3-rarefaction curve $\mathcal{S}_3(U_0)$ consisting of all left-hand states $U$ that can be connected to $U_0$ by a rarefaction wave is
	\begin{equation*}
		\begin{aligned}
			\mathcal{R}_3(U_0):\ &\rho=\rho(p,U_0)=\rho_0\left( \frac{p}{p_0}\right) ^{\frac{1}{\gamma}},\ p\leq p_0\\
			&u=u(p,U_0)=u_0+\frac{2a_0}{\gamma-1}\left[ \left( \frac{p}{p_0}\right) ^{\frac{\gamma-1}{2\gamma}}-1\right] ,\ p\leq p_0\\
		\end{aligned}
	\end{equation*}
\end{Proposition}

The characteristic domain of $\lambda_2=u$ is linearly degenerate. The elementary wave associated with the $\lambda_2$ characteristic domain is the contact discontinuity.
\begin{Proposition}[contact discontinuity curves]
	Given a left-hand state $U_0$, the 2-contact discontinuity curve $\mathcal{R}_1(U_0)$ consisting of all right-hand states $U$ that can be connected to $U_0$ by a contact discontinuity is
	\begin{equation*}
		\mathcal{C}_2(U_0):\ u=u_0,\ p=p_0.
	\end{equation*}
\end{Proposition}

\section{Structures of Riemann solution}\label{Structures of Riemann solution}
\begin{figure}[htbp]\label{figure: general structure}
	\centering
	\subfigure{
		\includegraphics[width=0.45\linewidth]{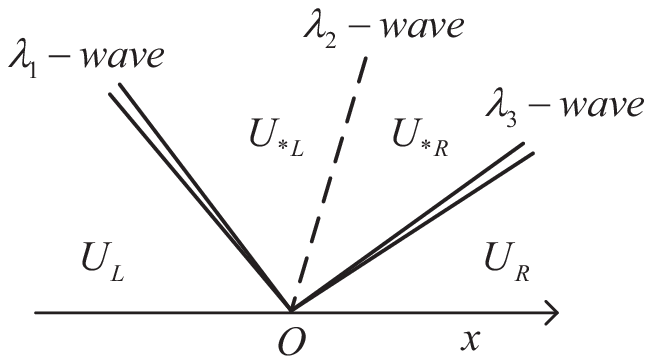}
	}
	\quad
	\subfigure{
		\includegraphics[width=0.47\linewidth]{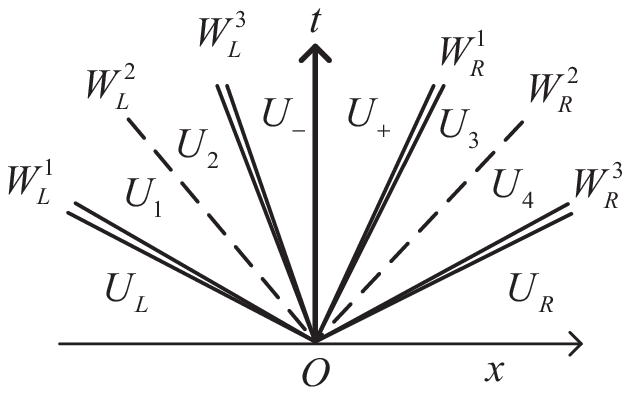}
	}
	\caption{The general structure of Riemann solutions. Left: associate Riemann problem (\ref{equ: Riemann problem 2}). Right: Riemann problem (\ref{equ: associate Riemann problem}).}
\end{figure}

The Riemann problem is a Cauchy problem of (\ref{equ: governing equation}) with piecewise constant initial conditions:
\begin{equation}\label{equ: Riemann problem 2}
	\left\lbrace 
	\begin{aligned}
		&\frac{\partial U}{\partial t}+\frac{\partial F}{\partial x}=\delta(x)S,\\
		&U(x,0)=\begin{cases}
			U_L,\ x<0\\
			U_R,\ x>0
		\end{cases}
	\end{aligned}
	\right. 
\end{equation}
The Riemann problem of classical Euler equations with the same initial value, which is also called associate Riemann problem, is
\begin{equation}\label{equ: associate Riemann problem}
	\left\lbrace 
	\begin{aligned}
		&\frac{\partial U}{\partial t}+\frac{\partial F}{\partial x}=0,\\
		&U(x,0)=\begin{cases}
			U_L,\ x<0\\
			U_R,\ x>0
		\end{cases}            
	\end{aligned}
	\right. 
\end{equation}\\
The solution of the associative Riemann problem (\ref{equ: associate Riemann problem}) is self-similar and contains four constant regions, which are divided by three classical elementary waves, as shown in the left figure of Figure\ref{figure: general structure}. Although the governing equation (\ref{equ: governing equation}) contains source terms, it is clearly still self-similar, thus we only consider self-similar solutions of the Riemann problem (\ref{equ: Riemann problem 2}). The following assumption will be needed throughout the paper.
\begin{Assumption}\label{assumption:self-similar}
	The solutions of Riemann problem (\ref{equ: Riemann problem 2}) is self-similar
	\begin{equation*}
		U(x,t)=U\left( \frac{x}{t}\right).
	\end{equation*}
\end{Assumption}
The self-similar solution of Riemann problem (\ref{equ: Riemann problem 2}) and associate Riemann problem (\ref{equ: associate Riemann problem}) with the initial value condition $U_L$ and $U_R$ are denoted as $U^R\left( \frac{x}{t},U_L,U_R\right),\ U^{A-R}\left( \frac{x}{t},U_L,U_R\right)$, respectively. Under Assumption\ref{assumption:self-similar}, $U_-=U^R\left( \frac{x}{t}=0-,U_L,U_R\right)$ and $U_+=U^R\left( \frac{x}{t}=0+,U_L,U_R\right)$ are both time-independent vectors. 

Double classical Riemann problem(CRP) framework will be used to analyze the structures of Riemann solutions in this section.
\begin{Proposition}[Double CRP framework]\label{proposition: double CRPs framework}
	The self-similar solution of the Riemann problem consists of seven discontinuities at most. They are a stationary wave at $x=0$, two genuinely nonlinear waves and a contact discontinuity left to $x=0$, two genuinely nonlinear waves and a contact discontinuity right to $x=0$ respectively.
\end{Proposition}
According to Proposition\ref{proposition: double CRPs framework}, the general structure of Riemann solutions is shown in the right figure of Figure\ref{figure: general structure}. The elementary waves on the left and right sides are denoted as $W_L^1$, $W_L^2$, $W_L^3$ and $W_R^1$, $W_R^2$, $W_R^3$, respectively. 
$W_L^1$, $W_L^3$, $W_R^1$ and $W_R^3$ can be shock waves or rarefaction waves, and both $W_L^2$ and $W_R^2$ are contact discontinuities. The t-axis is a stationary wave and $U_-\ne U_+$.
We follow the notation of \cite{lefloch2003riemann,thanh2009riemann} to express the structure of Riemann solutions. For example, $\mathcal{S}(U_L,U_1)$ and $\mathcal{R}(U_L,U_1)$ mean that two states $U_L$ and $U_1$ are connected by a shock and a rarefaction wave, respectively. The symbol $\mathcal{W}(U_L,U_1)$ means that  $U_L$ and $U_1$ are connected by a shock wave or a rarefaction wave.
$\mathcal{C}(U_1,U_2)$ means that $U_1$ and $U_2$ are connected by a contact discontinuity, and $\mathcal{D}(U_-,U_+)$ means that $U_-$ and $U_+$ are connected by a stationary wave. The connection between multiple elementary waves and constant regions is represented by the symbol "$\oplus $". For example, $\mathcal{S}(U_L,U_1)\oplus  \mathcal{C}(U_1,U_2)$ means that $U_L$ and $U_1$ are connected by a shock and a shock wave, followed by a contact discontinuity with the left-hand state $U_1$ and the right-hand state $U_2$.

\begin{Theorem}\label{theorem: structures}
	There are seven structures of the Riemann solution:
	\begin{itemize}
		\item [(i)]non-choked structures:
		\begin{itemize}
			\item[] Type1: $\mathcal{W}(U_L,U_-)\oplus \mathcal{D}(U_-,U_+)\oplus \mathcal{C}(U_+,U_4)\oplus \mathcal{W}(U_4,U_R)$.
			\item[] Type2: $\mathcal{D}(U_L,U_+)\oplus \mathcal{W}(U_+,U_3)\oplus \mathcal{C}(U_3,U_4)\oplus \mathcal{W}(U_4,U_R)$.
		\end{itemize}
		\item [(ii)]choked structures:
		\begin{itemize}
			\item[] Type3: $\mathcal{W}(U_L,U_-)\oplus \mathcal{D}(U_-,U_+)\oplus \mathcal{R}(U_+,U_3)\oplus \mathcal{C}(U_3,U_4)\oplus \mathcal{W}(U_4,U_R)$ with $M_+=1$.
			\item[] Type4: $\mathcal{D}(U_L,U_+)\oplus \mathcal{R}(U_+,U_3)\oplus \mathcal{C}(U_3,U_4)\oplus \mathcal{W}(U_4,U_R)$ with $M_+=1$.
			\item[] Type5: $\mathcal{R}(U_L,U_-)\oplus \mathcal{D}(U_-,U_+)\oplus \mathcal{W}(U_+,U_3)\oplus \mathcal{C}(U_3,U_4)\oplus \mathcal{W}(U_4,U_R)$ with $M_-=1$.
			\item[] Type6: $\mathcal{R}(U_L,U_-)\oplus \mathcal{D}(U_-,U_+)\oplus \mathcal{C}(U_+,U_4)\oplus \mathcal{W}(U_4,U_R)$ with $M_-=1$.
			\item[] Type7: $\mathcal{R}(U_L,U_-)\oplus \mathcal{D}(U_-,U_+)\oplus \mathcal{R}(U_+,U_3)\oplus \mathcal{C}(U_3,U_4)\oplus \mathcal{W}(U_4,U_R)$ with $M_-=M_+=1$.
		\end{itemize}
	\end{itemize}		
\end{Theorem}
\begin{Proof}\label{proof}
	Following the idea of double CRP framework, we first give the relations between the waves in the left and right half-planes with respect to states $U_-$ and $U_+$ respectively, and then combine these two sets of waves to form the structure of global Riemann solution under Criterion\ref{criterion1}.
	
	We begin by considering the structure of the left half-plane. 
	
	We use a contradiction to prove that the wave $W_L^3$ do not exist.
	If the assertion would not hold, then it follows from $u_-+a_->u_->0$ that $W_L^3$ cannot be a rarefaction wave, therefore $W_L^3$ is a shock wave. Let $\sigma_{L,3}$ denotes the speed of $W_L^3$. It follows from Lax entropy condition that $\sigma_{L,3}>u_-+a_->0$, which contradicts $\sigma_{L,3} <0$.
	We have proved that $W_L^3$ does not exist, then it follows from $u_->0$ that the wave $W_L^2$ does not exist either.
	
	We now consider the wave $W_L^1$ in the following three cases.
	
	First, if $0<M_-<1$, then it is clear that wave $W_L^1$ is possible in the solution.
	
	Second, if $M_->1$, then we assert that $W_L^1$ does not exist. Otherwise, if $W_L^1$ exists, then it must be a shock wave with negative speed. Let $\sigma_{L,1}$ denotes the speed of $W_L^1$. According to Lax entropy condition, we have $0>\sigma_{L,1}>u_--a_->0$, which is impossible. 
	
	Finally, if $M_-=1$, then we assert that $W_L^1$ is a rarefaction wave. If $W_L^1$ exists and is a shock wave, then Lax entropy condition shows that $\sigma_{L,1}>u_--a_-=0$, which is contrary to $\sigma_{L,1}\leq 0$. In addition, the nonexistence of $W_L^1$ can be regarded as the degeneration of a rarefaction wave, hence $W_L^1$ is a rarefaction wave. 
	
	We have thus related the structure of left half-plane to the value of $M_-$ as follows.
	\begin{itemize}
		\item [(1)] If $0<M_-<1$, then the Riemann solution contains only a classical nonlinear wave $W_L^1$ in the left half-plane.
		\item [(2)] If $M_-=1$, then the Riemann solution contains only a rarefaction wave $W_L^1$ in the left half-plane.
		\item [(3)] If $M_->1$, then the Riemann solutions do not contain waves in the left half-plane.
	\end{itemize}
	
	Using a similar approach, we can then prove the relations between elementary waves and $U_+$ in the right half-plane as follows.
	\begin{itemize}
		\item [(4)] If $0<M_+<1$, then the structure of Riemann solution in the right half-plane is $\mathcal{C}(U_+,U_4)\oplus \mathcal{W}(U_4,U_R)$.
		\item [(5)] If $M_+=1$, then the structure of Riemann solution in the right half-plane is $\mathcal{W}(U_+,U_3)\oplus \mathcal{C}(U_3,U_4)\oplus \mathcal{W}(U_4,U_R)$.
		\item [(6)] If $M_+>1$, then the structure of Riemann solution in the right half-plane is $\mathcal{R}(U_+,U_3)\oplus \mathcal{C}(U_3,U_4)\oplus \mathcal{W}(U_4,U_R)$.
	\end{itemize}
	
	The final step of the proof is to couple the left and right parts of elementary waves under Criterion\ref{criterion2}. According to Criterion\ref{criterion2}, there are two cases of $M_-$ and $M_+$ for non-choked structures. If $0<M_-,M_+<1$, the structure of solution is
	\begin{equation*}
		\mathcal{W}(U_L,U_-)\oplus \mathcal{D}(U_-,U_+)\oplus \mathcal{C}(U_+,U_4)\oplus \mathcal{W}(U_4,U_R).
	\end{equation*}
	If $M_-,M_+>1$, the structure of solution is
	\begin{equation*}
		\mathcal{D}(U_L,U_+)\oplus \mathcal{W}(U_+,U_3)\oplus \mathcal{C}(U_3,U_4)\oplus \mathcal{W}(U_4,U_R).
	\end{equation*}	
	The structure of choked solution is different for different signs of $k$. When $k>0$, it follows from \ref{table: left and right Mach numbers} that $M_-=M_1^*,M_+=1$ or $M_-=M_2^*,M_+=1$, then the two structures of choked solution are
	\begin{equation*}
		\begin{aligned}
			&\mathcal{W}(U_L,U_-)\oplus \mathcal{D}(U_-,U_+)\oplus \mathcal{R}(U_+,U_3)\oplus \mathcal{C}(U_3,U_4)\oplus \mathcal{W}(U_4,U_R),\\
			&\mathcal{D}(U_L,U_+)\oplus \mathcal{R}(U_+,U_3)\oplus \mathcal{C}(U_3,U_4)\oplus \mathcal{W}(U_4,U_R).
		\end{aligned}
	\end{equation*}
	When $k<0$, the choked solution satisfies $M_-=1,M_+=M_1^{**}$ or $M_-=1,M_+=M_2^{**}$, and its structures are
	\begin{equation*}
		\begin{aligned}
			&\mathcal{R}(U_L,U_-)\oplus \mathcal{D}(U_-,U_+)\oplus \mathcal{W}(U_+,U_3)\oplus \mathcal{C}(U_3,U_4)\oplus \mathcal{W}(U_4,U_R),\\
			&\mathcal{R}(U_L,U_-)\oplus \mathcal{D}(U_-,U_+)\oplus \mathcal{C}(U_+,U_4)\oplus \mathcal{W}(U_4,U_R).
		\end{aligned}
	\end{equation*}
	When $k=0$, the choked solution satisfies $M_-=M_+=1$, and its structure is
	\begin{equation*}
		\mathcal{R}(U_L,U_-)\oplus \mathcal{D}(U_-,U_+)\oplus \mathcal{R}(U_+,U_3)\oplus \mathcal{C}(U_3,U_4)\oplus \mathcal{W}(U_4,U_R).
	\end{equation*}
\end{Proof}

\begin{figure}[htbp]
	\label{figure: non-choked structures}
	\centering
	\subfigure[Type1]{
		\includegraphics[width=11em]{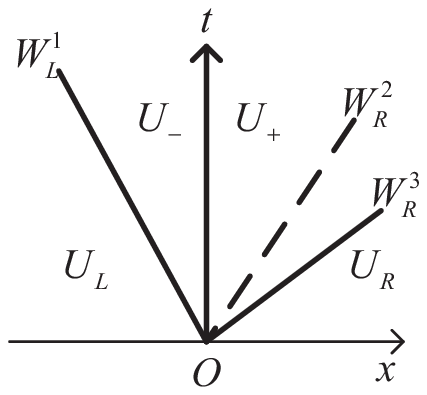}
	}
	\quad
	\subfigure[Type2]{ 
		\includegraphics[width=11em]{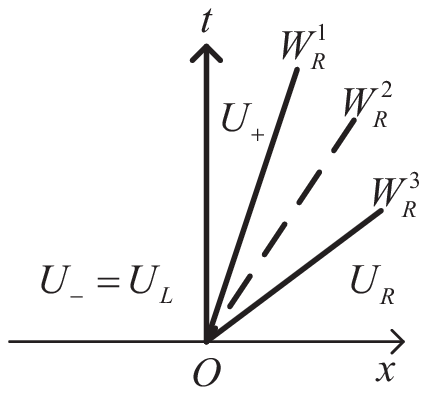}
	}
	\caption{The non-choked structures of Riemann solutions.}
\end{figure}
\begin{figure}[htbp]
	\label{figure: choked structures}
	\centering
	\subfigure[Type3]{
		\includegraphics[width=11em]{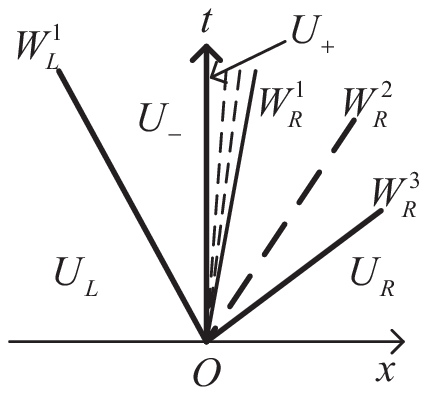}
	}
	\quad
	\subfigure[Type4]{
		\includegraphics[width=11em]{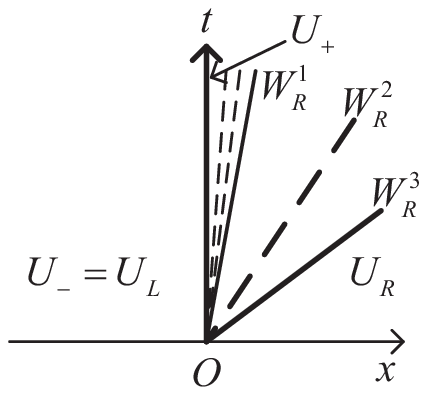}
	}
	\quad
	\subfigure[Type5]{
		\includegraphics[width=11em]{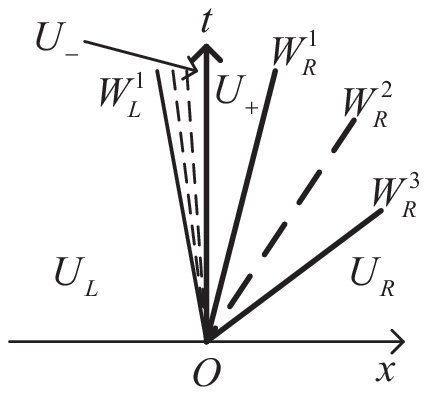}
	}
	\quad
	\subfigure[Type6]{
		\includegraphics[width=11em]{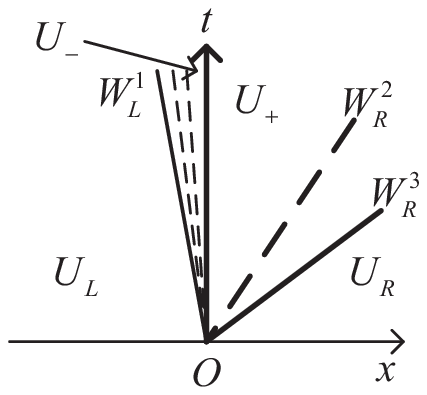}
	}
	\quad
	\subfigure[Type7]{
		\includegraphics[width=11em]{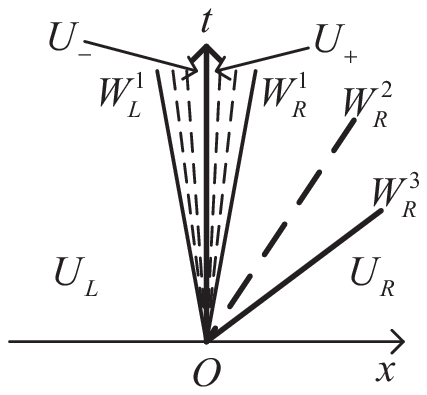}
	}
	\caption{The choked structures of Riemann solutions.}
\end{figure}

The sketches of these seven structures are shown in Figure\ref{figure: non-choked structures,figure: choked structures}. 
The non-choked structure is in the form of a stationary wave bisecting a certain constant region of a classical Riemann solution, and these two new regions satisfy the jump relation equation (\ref{equ: jump relation 1}).
The choked structure is more complex, where the elementary wave associated with a certain characteristic domain appears twice, and one of these two waves is a rarefaction wave and in contact with the stationary wave. 
From the proof of Thoerem\ref{theorem: structures} we have obtained all possible structures of Riemann solutions for the three ranges of coefficients $k$, as follows.
\begin{Corollary}\mbox{}
	\begin{itemize}
		\item[(i)] If $k>0$, then all possible structures of the Riemann solutions are Type1, Type2, Type3 and Type4.
		\item[(ii)] If $k<0$, then all possible structures of the Riemann solutions are Type1, Type2, Type5 and Type6.
		\item[(iii)] If $k=0$, then all possible structures of the Riemann solutions are Type1, Type2 and Type7.	 
	\end{itemize}
\end{Corollary}
The possible structures of solutions for different $k$ and their corresponding ranges of Mach numbers on either side of stationary wave are shown in table\ref{table: Mach number}.

\begin{table}[htbp]\label{table: Mach number}
	\centering
	\caption{For each $k$, the possible structure of the Riemann solution and the range of Mach numbers on both sides of the stationary wave.}
	\begin{tabular}{cccc}
		\toprule
		Coefficients & Solution structure & Upstream Mach number &Downstream Mach number\\
		\midrule
		\specialrule{0em}{3pt}{3pt}
		\multirow{4}{*}{$k>0$}
		&Type1& $0<M_-<M^*_1$&$0<M_+<1$\\	
		\specialrule{0em}{1pt}{1pt}	
		&Type2& $M_->M^*_2$&$1<M_+<M^*_3$\\		
		\specialrule{0em}{1pt}{1pt}	
		&Type3& $M_-=M^*_1$&$M_+=1$\\	
		\specialrule{0em}{1pt}{1pt}		
		&Type4& $M_-=M^*_2$&$M_+=1$\\
		\specialrule{0em}{3pt}{3pt}
		\midrule
		\specialrule{0em}{3pt}{3pt}
		\multirow{4}{*}{$k<0$}
		&Type1& $0<M_-<1$&$0<M_+<M^{**}_1$\\
		\specialrule{0em}{1pt}{1pt}			
		&Type2& $1<M_-<M^{**}_3$&$M_+>M^{**}_2$\\
		\specialrule{0em}{1pt}{1pt}			
		&Type5& $M_-=1$&$M_+=M^{**}_1$\\
		\specialrule{0em}{1pt}{1pt}			
		&Type6& $M_-=1$&$M_+=M^{**}_2$\\
		\specialrule{0em}{3pt}{3pt}
		\midrule
		\specialrule{0em}{3pt}{3pt}
		\multirow{3}{*}{$k=0$}
		&Type1& $0<M_-<1$&$0<M_+<1$\\
		\specialrule{0em}{1pt}{1pt}		
		&Type2& $M_->1$&$M_+>1$\\	
		\specialrule{0em}{1pt}{1pt}		
		&Type7& $M_-=1$&$M_+=1$\\
		\specialrule{0em}{3pt}{3pt}
		\bottomrule
	\end{tabular}
\end{table}

In some special cases, there are two admissible branches of the wave curve of the stationary wave. We show below that it does not lead to the Riemann problem being ill-posed.

The case of multiple solutions of the wave curve \ref{proposition: stationary wave curve1} satisfies
\begin{equation}
	k>0,M_-=M_1^*\text{ or }M_1^{**},M_+=1.
\end{equation}

\begin{Property}\label{property: uniquenesss1}
	Suppose the initial values of the Riemann solution $U_1(x,t)$ are $U_{L,1}$ and $U_{R,1}$ and $M_-=M_1^*,M_+=1$ are satisfied, and the initial values of the Riemann solution $U_2(x,t)$ are $U_{L,2}$ and $U_{R,2}$ and $M_-=M_1^{**},M_+=1$ are satisfied. If $U_{L,1}=U_{L,2},U_{R,1}=U_{R,2}$, then \[U_1(x,t)=U_2(x,t),\forall x\in \mathbb{R},t\geq 0.\]
\end{Property}
\begin{Proof}
	The structure of $U_1$ is Type3, and the structure of $U_2$ is Type4.
	If $U_{L,1}=U_{L,2}$, then $M_{L,1}=M_1^{**}$.
	
	For $U_1$, we have $M_-=M_1^*$.
	From the fact that $M_1^*$ and $M_1^{**}$ satisfy the Prandtl relation
	\[
	(M_1^{**})^2=\frac{1+[(\gamma-1)/2(M_1^*)^2]}{\gamma (M_1^*)^2-(\gamma-1)/2}
	\]
	we know that the wave $W_L^1$ is a normal shock wave, therefore
	\[
	U_1(x,t)=U_2(x,t)\equiv U_L, \forall x<0,t\geq 0.
	\]
	By $U_1(0+,t)=U_2(0+,t)$ and the uniqueness of classical Riemann solution, we have
	\[
	U_1(x,t)=U_2(x,t), \forall x>0,t\geq 0.
	\]
\end{Proof}

The case of multiple solutions of the wave curve Proposition\ref{proposition: stationary wave curve2} satisfies
\begin{equation}
	k<0,M_-=1,M_+=M_2^*\text{ or }M_2^{**}.
\end{equation}

\begin{Property}\label{property: uniquenesss2}
	Suppose the initial values of the Riemann solution $U_1(x,t)$ are $U_{L,1}$ and $U_{R,1}$ and $M_-=1,M_+=M_2^{*}$ are satisfied, and the initial values of the Riemann solution $U_2(x,t)$ are $U_{L,2}$ and $U_{R,2}$ and $M_-=1,M_+=M_2^{**}$ are satisfied. If $U_{L,1}=U_{L,2},U_{R,1}=U_{R,2}$, then \[U_1(x,t)=U_2(x,t),\forall x\in \mathbb{R},t\geq 0.\]
\end{Property}
The proof of this result is quite similar to that of Property\ref{property: uniquenesss1} and so is omitted.
Property\ref{property: uniquenesss1} and Property\ref{property: uniquenesss2} show that the multiple admissible branches of stationary wave curves do not lead to the multiple solutions of the Riemann problem.

\section{Conclusions}\label{Conclusions}
The Riemann problem of Euler equations with a singular source was studied.
We proposed an eigenvalue-based monotonicity criterion for selecting the physical solution of stationary waves.
The criterion is also valid for singular source-induced discontinuities in other equations, such as the shallow water equations. 
We analyzed all possible structures of Riemann solution under the double CRP framework.
The theoretical conclusions in this paper on the stationary wave and the Riemann problem are very helpful in the design of numerical schemes, and the related work will be given in a subsequent paper.

\section*{Acknowledgments}
This work was supported by the National Natural Science Foundation of China (Grant No.12101029) and Postdoctoral Science Foundation of China (Grant No.2020M680283).

\bibliographystyle{plain}
\bibliography{references}

\end{document}